\documentclass{article}

\usepackage[utf8]{inputenc}
\usepackage[T1]{fontenc}
\usepackage{lmodern}
\usepackage{graphicx}
\usepackage{color}
\usepackage{hyperref}
\usepackage{amsmath}
\usepackage{amsfonts}
\usepackage{epstopdf}
\usepackage[table]{xcolor}
\usepackage{float}
\usepackage{matlab}

\begin{document}
\hypersetup{%
linkcolor=black,
citecolor=black}

\begin{titlepage}
   \begin{center}
       \vspace*{1cm}

       \textbf{Dynamical Systems and Markov Chains}

       \vspace{0.5cm}
            
       \vspace{1.5cm}

       \textbf{Ricardo Frumento}

       \vfill

       \vspace{0.8cm}

       University of South Florida\\
       Tampa, FL\\
       Fall 2020
            
   \end{center}
\end{titlepage}

\tableofcontents

\section{Abstract}
\begin{par}
This project is going to work with one example of stochastic matrix to understand how Markov chains evolve and how to use them to make faster and better decisions only looking to the present state of the system.
\end{par}
\section{Motivation}
\begin{par}
Using the knowledge of eigenvalues and eigenvectors it is possible to study how dynamical systems described by difference equations like $$\textbf{x}_{k+1}=A\textbf{x}_k$$ evolve through time.\cite{LAlgebra} The applications can range from weather, ecological problems, computer science and physics. The purpose of this computation is to find a steady state or long-term behaviour for the system.
\end{par}
\subsection{Stochastic Systems}
\begin{par}
A stochastic process is a family of randomly evolving variables. \cite{Stochastic} This is an important subject because real world problems have the same characteristic. Bacterial growth, electrical system fluctuations and the movement of gas are all examples of random processes. There are several statistical tools to study this kind of progression, but this project will focus on Markov chains and how they describe the changes happening in the system. Markov chains have applications not only in pure sciences like physics, chemistry, and biology but also in applied sciences like solar irradiance, internet, and speech recognition.
\end{par}
\subsection{Markov Chains}
\begin{par}
One way to describe a system that the probability of each new event only depends on the previous event is to use a stochastic model known as Markov Chain. This predictions are attractive because it is not necessary to look at the past to make them, which have equivalent accuracy as other methods that require knowledge of previous events. \cite{Chains} But to use it, there are a few assumptions that need to be made before using the Markov Chains.
\end{par}
\section{Description and Solution}
\subsection{Initial Assumptions}
\begin{par}
To use the Markov Chain, the following circumstances need to be present
$$\begin{array}{l}
    \text{\textbullet Matrix A is diagonalisable}\\
    \text{\textbullet A has n linearly independent eigenvectors}\  \textbf{v}_1,\dots,\textbf{v}_n
\end{array}$$
Once this is true, the eigenvectors form a basis for $R^n$ it is possible to write any initial vector as $$\textbf{x}_0=c_1\textbf{v}_1+\dots+c_n\textbf{v}_n$$ As the system follows the difference equation $$\textbf{x}_{n+1}=A\textbf{x}_n$$ The first two therms are $$\textbf{x}_1=A\textbf{x}_0=c1A\textbf{v}_1+\dots+c_nA\textbf{v}_n$$ $$\textbf{x}_1=c_1\lambda_1\textbf{v}_1+\dots+c_n\lambda_n\textbf{v}_n$$ $$\textbf{x}_2=A\textbf{x}_1=c1\lambda_1A\textbf{v}_1+\dots+c_n\lambda_nA\textbf{v}_n$$ $$\textbf{x}_2=c_1\lambda_1^2\textbf{v}_1+\dots+c_n\lambda_n^2\textbf{v}_n$$ Generalizing $$\textbf{x}_k=c_1{\lambda_1}^k\textbf{v}_1+\dots+c_n{\lambda_n}^k\textbf{v}_n$$ Now that there is an expression to the $k^{th}$ state it is possible to study the long-therm behaviour, it is just a matter of learning what happens when $k\rightarrow\infty$. An example helps with the visualization
\end{par}
\subsection{Solving a System}
\begin{par}
Consider the following stochastic matrix and dynamical system $$\textbf{x}_{k+1}=A\textbf{x}_k$$ $$A=\left[\begin{array}{cc} .8 & .1 \\ .2 & .9 \end{array}\right]$$ Solving $$\det(A-\lambda I)=0$$ The eigenvalues are $1$ and $0.7$. Now for the eigenvectors $$(A-\lambda I)\textbf{v}=0$$ For $\lambda_1=1$ and $\lambda_2=0.7$ the eigenvectors are $$\textbf{v}_1=\left[\begin{array}{c} 1 \\ 2 \end{array}\right],\quad\textbf{v}_2=\left[\begin{array}{c} -1 \\ 1 \end{array}\right]$$ Now, using the expression found for the $k^{th}$ term of a dynamical system the steady state can be studied $$\textbf{x}_k=c_1{\lambda_1}^k\textbf{v}_1+\dots+c_n{\lambda_n}^k\textbf{v}_n=c_11^k\textbf{v}_1+\dots+c_n0.7^k\textbf{v}_n$$ When $k\rightarrow\infty$ as $|\lambda_2|<1$ $$\textbf{x}_k=c_1\textbf{v}_1$$ This shows that the system, when sufficient computations are made, falls on the eigenspace of $\textbf{v}_1$ no matter the initial state. To better visualize this, a few examples are worked using MatLab, the computations can be found in the Appendix, and the figure below shows all states after four determined initial states and two random. $$\textbf{p}=\left[\begin{array}{c} 2 \\ 4 \end{array}\right]$$ $$\textbf{q}=\left[\begin{array}{c} -6 \\ 6 \end{array}\right]$$ $$\textbf{r}=\left[\begin{array}{c} 7 \\ 2 \end{array}\right]$$ $$\textbf{s}=\left[\begin{array}{c} -5 \\ -4 \end{array}\right]$$ 
\end{par}
\begin{figure}[htp]
    \centering
    \includegraphics[width=5cm]{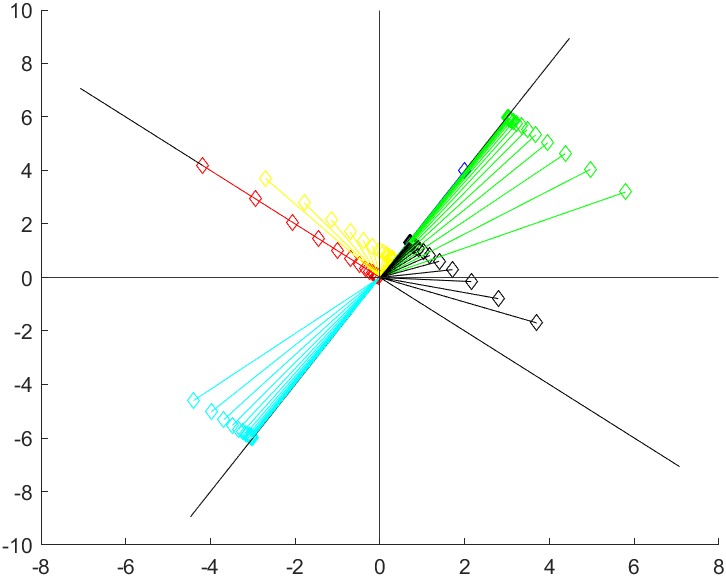}
    \caption{Evolution of the same system with different initial state}
    \label{fig:evolution}
\end{figure}
\begin{par}
In Figure 1, \textbf{p}, \textbf{q}, \textbf{r}, \textbf{s} are represented by the blue, red, green and cyan lines, respectively. The two randoms are depicted by the yellow and the black lines. The process to obtain this evolution is described for \textbf{p}, \textbf{q}, \textbf{r}, and \textbf{s} below.\\ \textbf{p} as initial state: $$\textbf{x}_0=\left[\begin{array}{c} 2 \\ 4 \end{array}\right]$$ $$\textbf{x}_1=A\textbf{x}_0=\left[\begin{array}{c} 2 \\ 4 \end{array}\right]$$ $$\vdots$$ $$\textbf{x}_{15}=A\textbf{x}_{14}=\left[\begin{array}{c} 2 \\ 4 \end{array}\right]$$ \textbf{q} as initial state: $$\textbf{x}_0=\left[\begin{array}{c} -6 \\ 6 \end{array}\right]$$ $$\textbf{x}_1=A\textbf{x}_0=\left[\begin{array}{c} -4.2 \\ 4.2 \end{array}\right]$$ $$\vdots$$ $$\textbf{x}_{15}=A\textbf{x}_{14}=\left[\begin{array}{c} -0.0285 \\ 0.0285 \end{array}\right]$$ \textbf{r} as initial state:  $$\textbf{x}_0=\left[\begin{array}{c} 7 \\ 2 \end{array}\right]$$ $$\textbf{x}_1=A\textbf{x}_0=\left[\begin{array}{c} 5.8 \\ 3.2 \end{array}\right]$$ $$\vdots$$ $$\textbf{x}_{15}=A\textbf{x}_{14}=\left[\begin{array}{c} 3.019 \\ 5.981 \end{array}\right]$$ \textbf{s} as initial state:  $$\textbf{x}_0=\left[\begin{array}{c} -5 \\ -4 \end{array}\right]$$ $$\textbf{x}_1=A\textbf{x}_0=\left[\begin{array}{c} -4.4 \\ -4.6 \end{array}\right]$$ $$\vdots$$ $$\textbf{x}_{15}=A\textbf{x}_{14}=\left[\begin{array}{c} -3.0095 \\ -5.9905 \end{array}\right]$$
\end{par}
\section{Discussion}
\begin{par}
Once the different initial state vectors are worked enough times, it is easy to assume they all have different steady state vectors, but after plotting the vectors and the eigenspaces of both eigenvectors the analysis is more complete and it is possible to say that all vectors are being transformed to a vector that is closer to the eigenspace of $\left[\begin{array}{c} 1 \\ 2 \end{array}\right]$ following a path that has the same direction as $\left[\begin{array}{c} -1 \\ 1 \end{array}\right]$. This is in accordance with what was expected once the equation $\textbf{x}_k=c_1\textbf{v}_1$ was found.
\end{par}
\begin{par}
Observing the characteristics of the stochastic matrix like the eigenvectors, it is important to note that they are orthogonal, meaning they form an orthogonal basis for $\Re^2$ and every initial state vector can be written as a linear combination of them. Therefore, once every initial state vector is written following the equation used to find $\textbf{x}_0$ and $\textbf{x}_k$ it is even more clear to where the long-term behavior is going, to the orthogonal projection of the initial state vector onto the eigenspace of the eigenvector associated with the eigenvalue $1$.
\end{par}
\section{Conclusion}
\begin{par}
Markov chains are a powerful tool for describing and analysing the evolution of dynamical systems. Once it is understood what transformations are happening on the vectors, it is fairly easy to analyse every state. For this project, it was only a matter of finding the orthogonal projection onto the eigenspace of one of the eigenvectors, but there are a lot of more interesting behaviours that can be studied. In this case, the stochastic matrix has an eigenvalue equal to one, and the other is bigger than zero and smaller than one.
\end{par}
\addcontentsline{toc}{section}{Bibliography}
\bibliographystyle{abbrv}
\bibliography{bibliographyMarkov.bib}

\begin{thebibliography}{1}

\bibitem{Chains}
Markov chains.
\newblock https://derrickchung.com/JAC/courses/Linear.105/.
\newblock Accessed on 2020-11-26.

\bibitem{Stochastic}
Stochastic process.
\newblock https://en.wikipedia.org/wiki/Stochastic\_process, Nov 2020.
\newblock Accessed on 2020-11-26.

\bibitem{LAlgebra}
J.~J.~M. David C.~Lay, Steven R.~Lay.
\newblock {\em Linear Algebra and Its Applications}.
\newblock Pearson Education, Inc., 2021.

\end{thebibliography}

\appendix
\section{MatLab Computations}
\begin{matlabcode}
A = [.8 .1; .2 .9]
\end{matlabcode}
\begin{matlaboutput}
A = 2x2    
    0.8000    0.1000
    0.2000    0.9000

\end{matlaboutput}
\begin{matlabcode}
[V,D] = eig(A)
\end{matlabcode}
\begin{matlaboutput}
V = 2x2    
   -0.7071   -0.4472
    0.7071   -0.8944

D = 2x2    
    0.7000         0
         0    1.0000

\end{matlaboutput}
\begin{matlabcode}
p = [2 4]';
q = [-6 6]';
r = [7 2]';
s = [-5 -4]';
t = randi([-10 10],1,2)';
u = randi([-10 10],1,2)';
xp = zeros(2,15);
xq = zeros(2,15);
xr = zeros(2,15);
xs = zeros(2,15);
xt = zeros(2,15);
xu = zeros(2,15);
time = zeros(1,15);
for i = 1:15
    time(i) = i;
    xp(:,i) = A^i*p;
    xq(:,i) = A^i*q;
    xr(:,i) = A^i*r;
    xs(:,i) = A^i*s;
    xt(:,i) = A^i*t;
    xu(:,i) = A^i*u;
end
\end{matlabcode}

\begin{matlabcode}
vec1 = V(:,1);
vec2 = V(:,2);
line1 = [-10*vec1 10*vec1];
line2 = [-10*vec2 10*vec2];
\end{matlabcode}

\begin{matlabcode}
clf
hold on
yline(0)
xline(0)
plot(line1(1,:),line1(2,:))
plot(line2(1,:),line2(2,:))
plot(xp(1,:), xp(2,:),'bd')
plot(xq(1,:), xq(2,:),'rd')
plot(xr(1,:), xr(2,:),'gd')
plot(xs(1,:), xs(2,:),'cd')
plot(xt(1,:), xt(2,:),'yd')
plot(xu(1,:), xu(2,:),'kd')
hold off
hold on
yline(0)
xline(0)
plot(line1(1,:),line1(2,:),'k')
plot(line2(1,:),line2(2,:),'k')
for j = 1:15
    vecp = [0 xp(1,j); 0 xp(2,j)];
    plot(vecp(1,:), vecp(2,:), 'b')
    vecq = [0 xq(1,j); 0 xq(2,j)];
    plot(vecq(1,:), vecq(2,:), 'r')
    vecr = [0 xr(1,j); 0 xr(2,j)];
    plot(vecr(1,:), vecr(2,:),'g')
    vecs = [0 xs(1,j); 0 xs(2,j)];
    plot(vecs(1,:), vecs(2,:), 'c')
    vect = [0 xt(1,j); 0 xt(2,j)];
    plot(vect(1,:), vect(2,:), 'y')
    vecu = [0 xu(1,j); 0 xu(2,j)];
    plot(vecu(1,:), vecu(2,:), 'k')
end
hold off
\end{matlabcode}
\end{document}